

\baselineskip=14pt
\parskip=10pt

\magnification=\magstephalf

\def\1{{\overline{1}}}
\def\2{{\overline{2}}}
\parindent=0pt
\overfullrule=0in

\def\frac#1#2{{#1 \over #2}}
\centerline
{\bf The $O(1/n^{85})$ Asymptotic expansion of OEIS sequence A85}
\bigskip
\centerline
{\it Shalosh B. EKHAD, Manuel KAUERS, and Doron ZEILBERGER}
\bigskip

\qquad \qquad {\it In honor of our two heroes, Don Knuth who was  $85$, a year ago, and Neil Sloane, who is now $85$}

On Oct. 10, 2024, {\bf exactly} $85$ years after his birth,
guru Neil Sloane gave a great Zoom talk: \hfill\break
{\tt https://vimeo.com/1018688642?share=copy} \quad . 

After the talk was over his many fans discussed the significance of the number $85$, and someone mentioned that the OEIS sequence A85 ({\tt https://oeis.org/A000085}) is a very
important sequence, the number of involutions of length $n$, and also, thanks to the famous RSK (Robinson-Schenstead-KNUTH) correspondence, the
number of Young tableaux with $n$ boxes. And indeed this sequence was important enough that guru Don Knuth, in his {\it bible}, The Art of Computer Programming
([Kn], pp. 62-64) spent three pages deriving, {\it by hand}, using {\it Laplace's method}, the asymptotic formula for the members of $A85$, first determined, using a different method
by Moser and Wyman ([MW]).  Namely, $t_n:=A85[n]$ (Eq. (53), p. 64, in [Kn]) is given by the following formula.
$$
t_n \, = \, \frac{1}{\sqrt{2}} n^{n/2}e^{-n/2+\sqrt{n}-1/4} \, \left ( 1+ \frac{7}{24} n^{-1/2} +O(n^{-3/4})  \right ) \quad .
$$

Knuth remarked:

\qquad {\it {\bf In principle}, the method we have used could be extended to any $O(n^{-k})$, for {\bf any} $k$,  instead of $O(n^{-3/4}$).}

In this modest tribute to our two heroes, we will do it  {\bf in practice} for the important number $k=85$.  Ideally one should be able to  use the Maple
package   \hfill\break
{\tt http://sites.math.rutgers.edu/\~{}zeilberg/tokhniot/AsyRec.txt} written by the third author, and executed by the
first author, that, alas, for some mysterious reason, can only go as far as $k=22$. Its extension, written in Sage by the second author, using the method in [KaJJ],
and executed by his computer, can go as far as one wishes, but in honor of the two gurus we will stop at $O(1/n^{85})$.

We have (to order $O(n^{-5})$)

{\bf Theorem}:
$$
t_n \, = \, \frac{1}{\sqrt{2}} n^{n/2}e^{-n/2+\sqrt{n}-1/4} \, \cdot
$$
$$
(1+ \frac{7}{24 \sqrt{n}}-\frac{119}{1152 n}-\frac{7933}{414720 n^{\frac{3}{2}}}+\frac{1967381}{39813120 n^{2}}-\frac{57200419}{1337720832 n^{\frac{5}{2}}}-\frac{562799}{47775744 n^{3}}-\frac{526420847}{40131624960 n^{\frac{7}{2}}}
$$
$$
+\frac{1856209}{573308928 n^{4}}-\frac{267645803}{2407897497600 n^{\frac{9}{2}}}\, + \,  O(\frac{1}{n^5}) ) \quad .
$$

To get the $O(\frac{1}{n^{85}})$  asymptotics, that is too long to be typyset in {\it humanese}, see the output file:

{\tt http://sites.math.rutgers.edu/\~{}zeilberg/tokhniot/A85asympt.txt} \quad .

{\bf Comments}

{\bf 1.} Our method is different from both the one used by Moser and Wyman  and the one used by Knuth (that he attributes to Laplace). We use the method exposited in [WZ], originally invented
by George D. Birkhoff and W. J. Trjitzinsky, and implemented by us in Maple and Sage.
We use the fact, mentioned in [Kn] (p. 62, Eq. (40) there), that  $t_n$ (i.e. $A85[n]$) satisfies the second-order recurrence
$$
t_n= t_{n-1}+ (n-1) t_{n-2} \quad .
$$

{\bf 2.} The combinatorial proof given in [Kn] (p. 62),   goes like this:

{\it 
A permutation is its own inverse if and only if its cycle form consists of one-cycles (fixed points) and two cycles (transpositions).
Since $t_{n-1}$ of the $t_n$ involutions have $(n)$ as a one-cycle, and since $t_{n-2}$ of them have $(jn)$ as a two-cyle, for fixed $j<n$, we obtain this formula.
}

We have two additional proofs, for what there are worth. Eq, $(41)$ of ([Kn], p.62) says:
$$
t_n \,= \, \sum_{k=1}^{ \lfloor n/2 \rfloor} \frac{n!}{(n-2k)!2^k k!} \quad,
$$

(as explained in [Kn], suppose there are $k$ two-cycles and $n-2k$ one-cycles. There are ${{n} \choose {2k}}$ ways to choose the fixed points,
and the multinomial coefficient $(2k)!/(2!)^k$ is the number of ways to arrange the other elements into $k$ distinguished transpositions;
dividing by $k!$ to make the transpositions indistinguishable.)

Now download the Maple package:

{\tt https://sites.math.rutgers.edu/\~{}zeilberg/tokhniot/EKHAD.txt} \quad,

and type

{\tt zeil(n!/((n-2*k)!*2**k*k!),k,n,N);} \quad,

and you would get the above recurrence, followed by its {\it proof certificate}.

Eq. $(42)$ in [Kn], p. 62,
$$
\sum_{n} t_nz^n/n! = e^{z+z^2/2} \quad ,
$$
yields yet-another proof. In the above-mentioned Maple package, type

{\tt AZd(n!*exp(z+z**2/2)/z**(n+1),z,n,N);} \quad,

and you would get the recurrence,  followed by its {\it proof certificate}.

{\bf 3.} This article is accompanied by a Maple package {\tt A85.txt} available from

{\tt https://sites.math.rutgers.edu/\~{}zeilberg/tokhniot/A85.txt} \quad,

that lets you find the asymptotic expansion to order $O(1/n^{(k+1)/2})$ for {\it any} $1 \leq k \leq 169$. Just type

`{\tt AsyI(n,k);}' if $n$ is numeric, and `{\tt AsyIs(n,k);}', if $n$ is symbolic.

Procedures {\tt SeqWn(N)} outputs the list of first $N$ terms of OEIS sequence {\tt A85}.

Just for fun, to get $A85[1000]$, (the  number of involutions of length $1000$) type:

{\tt SeqWn(1000)[1000];} \quad ,

getting, {\it right away},  the number in the following output file:

{\tt https://sites.math.rutgers.edu/\~{}zeilberg/tokhniot/A85w1000.txt} \quad.

This is a 1296-digit number, that in floating point to twenty digits is: 

$2.1439289538422655419 \cdot 10^{1296} \quad$ .

The approximation from using $O(1/n)$ asymptotics is gotten by typing

{\tt evalf(AsyI(1000,1),20);}

yielding
$$
2.1441496003431008422 \cdot 10^{1296} \quad .
$$
The ratio to the exact value is:
$$
1.0001029168902448312
$$

On the other hand with {\tt evalf(AsyI(1000,30),20);} (the $O(1/n^{31/2})$) formula), the ratio is:
$$
0.9999999999999999999999999996\dots \quad .
$$

{\it Pas mal!} \quad.

{\bf References}

[Kn] Donald E. Knuth, {\it ``The Art of Computer Programming volume 3: Sorting and Searching''; Second Edition}, Addison-Wesley, 1998.

[KaJJ] Manuel Kauers, Maximilian Jaroschek, Fredrik Johansson, 
{\it Ore polynomials in Sage}, in  ``Computer Algebra and Polynomials'', Springer LNCS {\bf 8942}, pages 105-125, 2015.

[MW] Leo Moser and Max Wyman, {\it On solutions of $x^d=1$ in symmetric groups}, Canadian Journal of Mathematics {\bf 7} (1955),
159-168.

[WZ] Jet Wimp and Doron Zeilberger, {\it Resurrecting the Asymptotics of Linear Recurrences}, J. of Math. Anal. Appl. {\bf 111} (1985), 162-176.

\bigskip
\hrule
\bigskip
Shalosh B. Ekhad, c/o D. Zeilberger, Department of Mathematics, Rutgers University (New Brunswick), Hill Center-Busch Campus, 110 Frelinghuysen
Rd., Piscataway, NJ 08854-8019, USA. \hfill\break
Email: {\tt ShaloshBEkhad at gmail dot com}   \quad .
\bigskip
Manuel Kauers, Institute for Algebra, J. Kepler University Linz, Austria \hfill\break
Email: {\tt manuel dot kauers at jku dot at}
\bigskip
Doron Zeilberger, Department of Mathematics, Rutgers University (New Brunswick), Hill Center-Busch Campus, 110 Frelinghuysen
Rd., Piscataway, NJ 08854-8019, USA. \hfill\break
Email: {\tt DoronZeil at gmail  dot com}   \quad .

\bigskip
\hrule
\bigskip

Oct. 20, 2024.

\end